\documentclass[11pt]{amsart}

\usepackage[T1]{fontenc}
\usepackage{lmodern}
\usepackage{microtype}
\usepackage{amsmath,amssymb}
\usepackage[a4paper, margin=3cm]{geometry}
\usepackage[colorlinks=true,linkcolor=blue,citecolor=blue,urlcolor=blue]{hyperref}

\newtheorem{theorem}{Theorem}[section]
\newtheorem{lemma}[theorem]{Lemma}

\newtheorem{conj}[theorem]{Conjecture}
\theoremstyle{remark}

\DeclareMathOperator{\conv}{conv}
\DeclareMathOperator{\relint}{relint}
\DeclareMathOperator{\vol}{vol}
\DeclareMathOperator{\width}{width}
\newcommand{\kslice}{S^{(k)}}
\newcommand{\knslice}{S_n^{(k)}}
\newcommand{\Z}{\mathbb Z}
\newcommand{\R}{\mathbb R}
\newcommand{\hstar}{h^*}
\newcommand{\lstar}{\ell^*}
\newcommand{\cA}{\mathcal A}
\newcommand{\eul}[2]{\genfrac{\langle}{\rangle}{0pt}{}{#1}{#2}}

\title[Wide lattice polytopes]{Lattice polytopes of large width have real-rooted
Ehrhart \(h^*\)-polynomials}
\author{Benjamin Nill}
\address[B.~Nill]{Faculty of Mathematics, Otto-von-Guericke-Universit\"at Magdeburg, Universit\"atsplatz 2, 39106 Magdeburg, Germany.}
\email{benjamin.nill@ovgu.de}


\begin{document}

\begin{abstract}
In this note we prove that in fixed dimension the Ehrhart \(h^*\)-polynomial
of a lattice polytope of sufficiently large lattice width is real-rooted. In particular, this implies strict log-concavity and
unimodality of the $h^*$-vector and answers a question of Averkov, Hofscheier
and Nill. For a lattice simplex we prove the analogous statement
for its local \(h^*\)-polynomial, also called box polynomial. The proofs were found using ChatGPT 5.6 Sol and follow essentially directly from a result by Basu and Oertel that for large enough lattice width counting lattice points approximates the volume.
\end{abstract}

\maketitle

\section{Introduction}

Let \(P\subset \R^d\) be a \(d\)-dimensional lattice polytope (i.e., its vertices are in $\Z^d$). Its Ehrhart series is
\[
 \sum_{n\geq0}|nP\cap \Z^d|t^n
 =\frac{\hstar(P;t)}{(1-t)^{d+1}},
 \qquad
 \hstar(P;t)=\sum_{j=0}^d h_j^*(P)t^j.
\]
The coefficients \(h_j^*(P)\) are nonnegative integers
\cite{Ehrhart1962,Stanley1980} and the properties of the coefficient vector, called the {\em $h^*$-vector} of $P$, have been intensively studied in the literature. We refer to \cite{Ferroni} for a survey of interesting examples and open questions.

For a fixed lattice polytope, it has been shown that the \(h^*\)-polynomials of its dilates eventually have distinct negative real roots, see \cite{Jochemko2018} as well as \cite{BrentiWelker2009, BeckStapledon2010, Higashitani2019}. In this note we show that this statement can be generalized to any family of lattice polytopes whose lattice widths go to infinity. For this, let us recall the definition of the {\em lattice width} of a convex body \(K\subset \R^d\):
\begin{equation}
 \width(K):=\min_{0\neq u\in (\Z^d)^*}
 \left(\max_{x\in K}u(x)-\min_{x\in K}u(x)\right).
 \label{eq:width}
\end{equation}

\begin{theorem}\label{thm:h}
For every \(d\geq1\) there exists \(C_d>0\) such that the
\(h^*\)-polynomial of every \(d\)-dimensional lattice polytope of
lattice width greater than \(C_d\) has positive coefficients and \(d\) distinct negative real
roots.
\end{theorem}

In particular, the $h^*$-vectors of lattice polytopes of large width are strictly log-concave and unimodal. This answers negatively Question~5 by Averkov, Hofscheier and the author in \cite{AverkovHofscheierNill2023}. Let us remark that in \cite{AverkovHofscheierNill2023} examples were given that show that large lattice width does not necessarily imply the lattice polytope to be IDP (i.e., to have the Integer Decomposition Property). We refer to \cite{Braun-Survey, Ferroni} for these notions.

We also show the same result for the local \(h^*\)-polynomial of a lattice simplex which for a simplex is often known under the name {\em box polynomial}. It is the enumerator polynomial of its open fundamental
parallelepiped, see \cite{BetkeMcMullen1985,Stanley1992,KatzStapledon2016, Solus2019,GustafssonSolus2020}.

\begin{theorem}\label{thm:box}
For every \(d\geq1\) there exists \(C_d^\Box>0\) such that the local
\(h^*\)-polynomial of every \(d\)-dimensional lattice simplex $S$ of
lattice width greater than \(C_d^\Box\) has positive coefficients $\ell^*_1(S), \ldots, \ell^*_d(S)$, a simple root at \(0\)
and \(d-1\) distinct negative real roots.
\end{theorem}

Also here, we get as consequences strict log-concavity and unimodality for its coefficient vector.

The proofs of both results are very short and follow essentially directly from a result by Basu and Oertel. 

Let us note that it is possible to make the constants $C_d$ and $C_d^\Box$ explicit, but the orders will be very large and very likely far off from the optimal ones. 

\begin{conj}
    Theorem~\ref{thm:box} also holds for lattice polytopes.
\end{conj}

As in contrast to lattice simplices, no analogous combinatorial counting interpretation for local $h^*$-polynomials of arbitrary lattice polytopes is known, the methods here do not directly generalize.

\subsection*{Acknowledgments}

The proofs were found using ChatGPT 5.6 Sol which also produced a first draft of this paper. The author is solely responsible for this final version. This work is funded by the Deutsche Forschungsgemeinschaft (DFG, German Research Foundation) – 539867500 as part of the research priority program Combinatorial Synergies. 

\section{Prerequisites}

\subsection{The relation between lattice point counting and the volume}

In Lemma~4.1.11 of the thesis of Timm Oertel \cite{Oertel2014} a beautiful inequality was elegantly proven that gives an explicit bound on the quotient $\frac{|K \cap \Z^d|}{\vol(K)}$ of a full-dimensional convex body $K$ containing the unimodular copy of a dilate of a standard cube. From this, an explicit approximation of $\frac{|K \cap \Z^d|}{\vol(K)}$ in terms of its lattice width can be deduced showing that this quotient converges to~$1$ if the lattice width becomes sufficiently large. This was published in the more general mixed-integer setting by Basu and Oertel in \cite[Lemma~3.8]{BasuOertel2017}. We only need the pure integer version which we state here. Basu and Oertel write that this result is "well-known".

\begin{lemma}[Basu and Oertel]\label{lem:count}
There exists a constant $N_d$ such that if \(K\) is a convex body in $\R^d$ with $\width(K) > c N_d$ for some positive real number $c$ then
\[ e^{-1/c}\le \frac{|K\cap\mathbb Z^d|}{\operatorname{vol}(K)} \le e^{1/c} \]
\end{lemma}

Hence, if a family of convex bodies $K_n$ satisfies $\width(K_n) \longrightarrow \infty$, then $\frac{|K_n\cap\mathbb Z^d|}{\operatorname{vol}(K_n)} \longrightarrow 1$, and -- as is easy to see -- also  $\frac{|\relint(K_n)\cap\mathbb Z^d|}{\operatorname{vol}(K_n)} \longrightarrow 1$.

\subsection{Approximating the Eulerian polynomial}

In this paper, we use the following shifted version of the {\em Eulerian polynomial}
\[
 \cA_d(t):=\sum_{j=1}^d\eul d{j-1}t^j
\]
where $\eul d {j-1}$ for $1\leq j\leq d$ is given as 
\[
 \eul d{j-1}:=
 \sum_{i=0}^j(-1)^{j-i}
 \binom{d+1}{j-i}i^d
\]

The Eulerian polynomial has \(d\) simple real roots, one equal to \(0\) and all the others
negative; see e.g. \cite{FoataSchutzenberger1970}. 

We will also need the following fact (see e.g. \cite{Marden}).

\begin{lemma}\label{lem:real}
If a real polynomial has only simple real roots, then this is also true after any sufficiently small real coefficient perturbation.
\end{lemma}

The idea is to use Rouch\'e's theorem that there is exactly one zero in each disk of a system of pairwise disjoint disks centered at the simple real roots. Since these disks are invariant under complex conjugation, the corresponding zeros are real; being counted with multiplicity one, they are simple.

\section{Proof of Theorem~\ref{thm:h}}

Let $d \ge 1$ be fixed. Suppose that \(P_n\) is a sequence of \(d\)-dimensional lattice
polytopes with \(\width(P_n)\to\infty\), and put
\(V_n:=\vol(P_n)\). Lemma~\ref{lem:count} gives for each \(1\leq i\leq d\) that 
\[
 \frac{|iP_n\cap \Z^d|}{\vol(i P_n)} \longrightarrow 1,
 \quad \text{ so }
 \frac{|iP_n\cap \Z^d|}{V_n}\longrightarrow i^d,
\]
when $n \longrightarrow \infty$. Note also that $V_n \longrightarrow \infty$ for $n \longrightarrow \infty$ (this is well-known, in our setting of lattice polytopes one could for instance use \cite{LagariasZiegler}), hence, with the convention $|0P_n \cap \Z^d| := 1$ we also have $\frac{|0P_n\cap \Z^d|}{V_n}\longrightarrow 0^d$ for $n \longrightarrow \infty$. The inversion formula for the Ehrhart series is for \(1\leq j\leq d\)
\[
 h_j^*(P_n)=
 \sum_{i=0}^j(-1)^{j-i}
 \binom{d+1}{j-i}|iP_n\cap \Z^d|.
\]

This shows the coefficientwise convergence
\[
 \frac{\hstar(P_n;t)}{V_n} = \frac{1}{V_n} + \sum_{j=1}^d \frac{h_j^*(P_n)}{V_n} t^j \longrightarrow \sum_{j=1}^d \left(\sum_{i=0}^j(-1)^{j-i}
 \binom{d+1}{j-i} i^d\right) t^j = \cA_d(t)
\]
for $n \longrightarrow \infty$ (recall that $h_0^*(P_n) = 1$).

By Lemma~\ref{lem:real} it follows that, for all sufficiently large \(n\) (where `large' depends by Lemma~\ref{lem:count} solely on $d$), 
\(\hstar(P_n;t)\) has \(d\) distinct real roots, one close to each root of \(\cA_d\), as well as positive coefficients. Since \(h_0^*(P_n)=1\) and all coefficients are
nonnegative, none of these roots is zero or positive.  They are therefore all negative.\hfill\qed

\section{Proof of Theorem~\ref{thm:box}}

Let \(S=\conv(v_0,\ldots,v_d)\subset \R^d\) be a lattice simplex and
let
\[
 D:=d!\vol(S)=
 \left|\det
 \begin{pmatrix}
 v_0&\cdots&v_d\\
 1&\cdots&1
 \end{pmatrix}\right|.
\]
The closed lattice parallelepiped spanned by $S \times {1}$ is
\[
 \Pi_S:=
 \left\{\sum_{i=0}^d\lambda_i(v_i,1):
 0\le\lambda_i\le1\right\}.
\]
Let us define for $1 \le k \le d$ the $k$-slice
\[\kslice := \Pi_S\cap (\R  ^d \times \{k\}) = \left\{\sum_{i=0}^d\lambda_i(v_i,1):
 0\le\lambda_i\le1,\; \lambda_0 + \cdots +\lambda_d = k\right\}.\]
It is well-known (e.g. \cite{Braun-Survey}) that the local $h^*$-polynomial of $S$
\[
 \lstar(S;t):=\sum_{k=1}^d \ell^*_k(S)t^k
\]
has coefficients
\[\ell^*_k(S):=|\relint(S^{(k)})\cap(\Z^d \times \{k\})|.\]

Let $k \in \{1, \ldots, d\}$. We let $V^{(k)}$ be the set containing all $k$-sums of the vertices of $S$. Note that $V^{(k)} \times \{k\}$ is the vertex set of $\kslice$. Let $u \in (\Z^d)^*$ be a non-zero integer functional. Without loss of generality let $u(v_0) \le u(v_1) \le \cdots \le u(v_d)$, hence, 
\[\max_{x\in S}u(x)-\min_{x\in S}u(x) = u(v_d)-u(v_0).\]
Note that as $u(v_{d-k+i})-u(v_i) \ge 0$ for $i=1,\ldots,k-1$, we have
\[\max_{x\in V^{(k)}}u(x)-\min_{x\in V^{(k)}} u(x)= (u(v_{d-k+1}) + \cdots + u(v_d)) - (u(v_0) + \cdots + u(v_{k-1})) \ge u(v_d)-u(v_0).\]
Hence, we get
\begin{equation}
\width\left(S^{(k)}\right) \ge \width(S).
\label{eq:in}
\end{equation}
 
Let $e_0, \ldots, e_d$ be the standard basis of $\R^{d+1}$. Hence, the linear map \(T:\R^{d+1}\to \R^{d+1}\) defined by \(T(e_i)=(v_i,1)\) for $i=0, \ldots, d$ has determinant $D$. 
Recall the hypersimplex 
\[
 H_k:=\left\{\lambda\in[0,1]^{d+1}:  \lambda_0+\cdots+\lambda_d=k\right\}.
\]
The volume of \(H_k\) with respect to the affine lattice given by the lattice points in its affine hull is $\frac{\eul d{k-1}}{d!}$ (this fact is attributed to Laplace in \cite{Li2012}). Note that $H_k$ gets mapped to $\kslice$ via $T$. As restricting $T$ to an affine map from the affine hull of $H_k$ to the affine hull of $\kslice$ also has determinant $D$, we get
\[\vol\left(\kslice\right)=\frac{D \eul d{k-1}}{d!}.\]

Now, let \(S_n\) be a sequence of \(d\)-dimensional lattice simplices
with \(\width(S_n)\to\infty\) for $n\longrightarrow \infty$, and put \(D_n:=d!\vol(S_n)\). As by \eqref{eq:in} for each \(1\leq k\leq d\) the widths of the $k$-slices also go to infinity, Lemma~\ref{lem:count} yields $\frac{|\relint(\knslice) \cap(\Z^d \times \{k\})|}{\operatorname{vol}(\knslice)} \longrightarrow 1$, and thus 
\[
 \frac{d!}{D_n}\ell^*_k(S_n)\longrightarrow\eul d{k-1}.
\]
Hence
\[
 \frac{d!}{D_n}\lstar(S_n;t)\longrightarrow\cA_d(t)
\]
coefficientwise.  By Lemma~\ref{lem:real} for all sufficiently large \(n\) (again, using Lemma~\ref{lem:count} the bound can be chosen to be dependent only on $d$) the box
polynomial therefore has one simple real root near each root of
\(\cA_d\). In particular, as the constant coefficient of the box polynomial vanishes, the root near
\(0\) is exactly \(0\), and all other roots are negative. Note that by this convergence, all coefficients $\ell^*_k(S_n)$ (for $k=1, \ldots, d)$ are eventually positive.\hfill\qed

\bibliographystyle{amsplain}
\bibliography{wide-hstar}

@article{AverkovHofscheierNill2023,
  author  = {Averkov, Gennadiy and Hofscheier, Johannes and Nill, Benjamin},
  title   = {Generalized flatness constants, spanning lattice polytopes, and the {G}romov width},
  journal = {Manuscripta Math.},
  volume  = {170},
  year    = {2023},
  pages   = {147--165},
  doi     = {10.1007/s00229-021-01363-x}
}

@article {LagariasZiegler,
    AUTHOR = {Lagarias, Jeffrey C. and Ziegler, G\"unter M.},
     TITLE = {Bounds for lattice polytopes containing a fixed number of
              interior points in a sublattice},
   JOURNAL = {Canad. J. Math.},
  FJOURNAL = {Canadian Journal of Mathematics. Journal Canadien de
              Math\'ematiques},
    VOLUME = {43},
      YEAR = {1991},
    NUMBER = {5},
     PAGES = {1022--1035},
      ISSN = {0008-414X,1496-4279},
   MRCLASS = {52C07 (11H06)},
  MRNUMBER = {1138580},
MRREVIEWER = {J.\ M.\ Wills},
       DOI = {10.4153/CJM-1991-058-4},
       URL = {https://doi.org/10.4153/CJM-1991-058-4},
}

@incollection {Braun-Survey,
    AUTHOR = {Braun, Benjamin},
     TITLE = {Unimodality problems in {E}hrhart theory},
 BOOKTITLE = {Recent trends in combinatorics},
    SERIES = {IMA Vol. Math. Appl.},
    VOLUME = {159},
     PAGES = {687--711},
 PUBLISHER = {Springer, [Cham]},
      YEAR = {2016},
doi={10.1007/978-3-319-24298-9_27}
}

@book {Marden,
    AUTHOR = {Marden, Morris},
     TITLE = {Geometry of polynomials},
    SERIES = {Mathematical Surveys},
    VOLUME = {No. 3},
   EDITION = {Second},
 PUBLISHER = {American Mathematical Society, Providence, RI},
      YEAR = {1966},
     PAGES = {xiii+243},
   MRCLASS = {30.11},
  MRNUMBER = {225972},
MRREVIEWER = {O.\ Shisha},
}

@article{Ferroni,
  title={Examples and counterexamples in {Ehrhart} theory},
  author={Ferroni, Luis and Higashitani, Akihiro},
  journal={preprint arXiv:2307.10852, to appear in EMS Surv. Math. Sci. },
  year={2024},
doi={10.48550/arXiv.2307.10852},
}

@article{BasuOertel2017,
  author  = {Basu, Amitabh and Oertel, Timm},
  title   = {Centerpoints: a link between optimization and convex geometry},
  journal = {SIAM J. Optim.},
  volume  = {27},
  year    = {2017},
  number  = {2},
  pages   = {866--889},
  doi     = {10.1137/16M1092908}
}

@article{BeckStapledon2010,
  author  = {Beck, Matthias and Stapledon, Alan},
  title   = {On the log-concavity of {H}ilbert series of {V}eronese subrings and {E}hrhart series},
  journal = {Math. Z.},
  volume  = {264},
  year    = {2010},
  number  = {1},
  pages   = {195--207},
  doi     = {10.1007/s00209-008-0458-7}
}

@article{BetkeMcMullen1985,
  author  = {Betke, Ulrich and McMullen, Peter},
  title   = {Lattice points in lattice polytopes},
  journal = {Monatsh. Math.},
  volume  = {99},
  year    = {1985},
  number  = {4},
  pages   = {253--265},
  doi     = {10.1007/BF01312545}
}

@article{BrentiWelker2009,
  author  = {Brenti, Francesco and Welker, Volkmar},
  title   = {The {V}eronese construction for formal power series and graded algebras},
  journal = {Adv. in Appl. Math.},
  volume  = {42},
  year    = {2009},
  number  = {4},
  pages   = {545--556},
  doi     = {10.1016/j.aam.2009.01.001}
}

@article{Ehrhart1962,
  author  = {Ehrhart, Eug{\`e}ne},
  title   = {Sur les poly{\`e}dres rationnels homoth{\'e}tiques {\`a} \(n\) dimensions},
  journal = {C. R. Acad. Sci. Paris},
  volume  = {254},
  year    = {1962},
  pages   = {616--618}
}

@book{FoataSchutzenberger1970,
  author    = {Foata, Dominique and Sch{\"u}tzenberger, Marcel-Paul},
  title     = {Th{\'e}orie g{\'e}om{\'e}trique des polyn{\^o}mes eul{\'e}riens},
  series    = {Lecture Notes in Mathematics},
  volume    = {138},
  publisher = {Springer-Verlag},
  address   = {Berlin-New York},
  year      = {1970},
  doi       = {10.1007/BFb0060799}
}

@article{GustafssonSolus2020,
  author  = {Gustafsson, Nils and Solus, Liam},
  title   = {Derangements, {E}hrhart theory, and local \(h\)-polynomials},
  journal = {Adv. Math.},
  volume  = {369},
  year    = {2020},
  pages   = {107169},
  doi     = {10.1016/j.aim.2020.107169}
}

@article{Higashitani2019,
  author  = {Higashitani, Akihiro},
  title   = {Unimodality of \(\delta\)-vectors of lattice polytopes and two related properties},
  journal = {European J. Math.},
  volume  = {5},
  year    = {2019},
  number  = {2},
  pages   = {333--355},
  doi     = {10.1007/s40879-018-0282-5}
}

@article{Jochemko2018,
  author  = {Jochemko, Katharina},
  title   = {On the real-rootedness of the {V}eronese construction for rational formal power series},
  journal = {Int. Math. Res. Not. IMRN},
  year    = {2018},
  number  = {15},
  pages   = {4780--4798},
  doi     = {10.1093/imrn/rnx027}
}

@article{KatzStapledon2016,
  author  = {Katz, Eric and Stapledon, Alan},
  title   = {Local \(h\)-polynomials, invariants of subdivisions, and mixed {E}hrhart theory},
  journal = {Adv. Math.},
  volume  = {286},
  year    = {2016},
  pages   = {181--239},
  doi     = {10.1016/j.aim.2015.09.010}
}

@article{Li2012,
  author  = {Li, Nan},
  title   = {Ehrhart \(h^*\)-vectors of hypersimplices},
  journal = {Discrete Comput. Geom.},
  volume  = {48},
  year    = {2012},
  number  = {4},
  pages   = {847--878},
  doi     = {10.1007/s00454-012-9452-2}
}

@phdthesis{Oertel2014,
  author       = {Oertel, Timm},
  title        = {Integer Convex Minimization in Low Dimensions},
  school       = {ETH Zurich},
  type         = {Doctoral Thesis},
  year         = {2014},
  doi          = {10.3929/ethz-a-010295887},
}

@incollection{Solus2019,
  author    = {Solus, Liam},
  title     = {Local \(h^*\)-polynomials of some weighted projective spaces},
  booktitle = {Algebraic and Geometric Combinatorics on Lattice Polytopes},
  publisher = {World Scientific},
  address   = {Hackensack, NJ},
  year      = {2019},
  pages     = {382--399}
}

@article{Stanley1980,
  author  = {Stanley, Richard P.},
  title   = {Decompositions of rational convex polytopes},
  journal = {Ann. Discrete Math.},
  volume  = {6},
  year    = {1980},
  pages   = {333--342},
  doi     = {10.1016/S0167-5060(08)70717-9}
}

@article{Stanley1992,
  author  = {Stanley, Richard P.},
  title   = {Subdivisions and local \(h\)-vectors},
  journal = {J. Amer. Math. Soc.},
  volume  = {5},
  year    = {1992},
  number  = {4},
  pages   = {805--851},
  doi     = {10.1090/S0894-0347-1992-1157293-9}
}

\end{document}